\documentclass[11pt,oneside]{amsart}
\usepackage{calc,ifthen,ulem, geometry,verbatim, color,graphicx, mathabx,mathrsfs,xspace,enumerate,amsmath,mathtools, amsthm, amssymb,yfonts}
\usepackage[usenames,dvipsnames,svgnames,table]{xcolor}
\usepackage[pdftex,bookmarks,colorlinks,breaklinks]{hyperref}  
\hypersetup{linkcolor=blue,citecolor=red,filecolor=dullmagenta,urlcolor=darkblue} 
\newtheorem{theorem}{Theorem}[section]
\newtheorem{corollary}[theorem]{Corollary}

\newtheorem{lemma}[theorem]{Lemma}
 
\newcommand{\cG}{\mathcal{G}}
\newcommand{\cH}{\mathcal{H}}
\newcommand{\cF}{\mathcal{F}}
\newcommand{\al}{\alpha}

\newcommand{\cc}{\mathcal}

\newcommand{\scr}{\mathscr}
\usepackage{mathtools}

\usepackage{enumitem}
\usepackage{bm}

  \newcommand\mult{\operatorname{\textup{{\fontfamily{ptm}\selectfont mult}}}}
  \newcommand\dg{\operatorname{\textup{{\fontfamily{ptm}\selectfont deg}}}}
  \newcommand\frc{\operatorname{\textup{{\fontfamily{ptm}\selectfont frac}}}}

     \makeatletter
      \def\@setcopyright{}
      \def\serieslogo@{}
      \makeatother
\begin{document}
   \author{Amin  Bahmanian}
   \address{Department of Mathematics,
  Illinois State University, Normal, IL USA 61790-4520}
\title[Connected Fair Detachments of Hypergraphs]{Connected Fair Detachments of Hypergraphs I}

   \begin{abstract}  
Let $\mathcal G$ be a hypergraph whose edges are colored.  An {\it $(\alpha,n)$-detachment} of $\mathcal G$ is a hypergraph obtained by splitting a vertex $\alpha$ into $n$ vertices, say $\alpha_1,\dots,\alpha_n$, and sharing the incident  edges  among the subvertices. A detachment is {\it fair} if the degree of vertices and multiplicity of edges are shared as evenly as possible among the subvertices within the whole hypergraph as well as within each color class. In this paper we solve an open problem from the 1970s by finding  necessary and sufficient conditions under which a $k$-edge-colored hypergraph $\cG$ has a fair detachment in which  each color class is connected.  Previously, this was not even known for the  case when $\cG$ is an arbitrary graph (i.e. 2-uniform hypergraph).  We exhibit the usefulness of our theorem by proving a variety of new results on  hypergraph decompositions (in Part I), and completing partial regular combinatorial structures (in Part II). 
   \end{abstract}
   \date{\today}
   \maketitle

\section{Introduction and  Statement of Main Result}
All hypergraphs under consideration are finite. We allow edges to occur multiple times, and we also allow each vertex to  have multiple occurrences within an edge. In other words, the edge set and every element of the edge set is a multiset. By an edge of the form $\{u_1^{m_1},\dots,u_s^{m_s}\}$, we mean an edge in which vertex $u_i$ occurs $m_i$ times for $1\leq i\leq s$.

Let $n$ be a positive integer. Given a hypergraph $\cH:=(V,E)$ and a vertex $\alpha$ of $\cH$, we obtain a hypergraph $\cF$ on $|V|+n-1$ vertices and $|E|$ edges in the following way.
\begin{enumerate} [label=(\roman*)]
\item Replace $\alpha$ in $\cG$ by $n$ new vertices $\alpha_1,\dots,\alpha_n$ in $\cF$;
\item Replace each edge $\{\alpha^p\}\cup U$ in $\cG$ by an edge $\{\alpha_1^{p_1},\dots,\alpha_n^{p_n}\}\cup U$ in $\cF$, where $p\geq 1, \alpha\notin U\subset V$ and $p_1+\dots+p_n=p$. 
\item Leave the remaining vertices and edges of $\cG$ intact and add them to $\cF$.
\end{enumerate} 
The resulting hypergraph $\cF$ is called an {\it $(\alpha,n)$-detachment}  of $\cG$, and $\cG$ is an {\it amalgamation} of $\cF$ (Amalgamation here is different from that of Ne\v{s}et\v{r}il \cite{MR556052}).  It is clear that $\cF$ is in
general far from being uniquely determined. 

There are two  active and somewhat independent lines of research on detachments of (hyper)graphs. In the first  line of research, which was initiated by Nash-Williams \cite{MR587625,MR810558,MR822773,MR916377,MR1104821}, an arbitrary graph is given, and the aim is to find conditions under which the detachment maintains certain structural properties such as edge-connectivity properties or that the detached vertices have prescribed degrees. This area has connections with submodular functions, combinatorial optimization and matroid theory \cite{MR1974484,MR2134419,MR1967880,MR2033306}.
For an overview of this line of work   together with  several other related results on variations of ``splitting off" operations, we refer the reader to  \cite{MR2223390}.

In the second line of research, which was initiated by Hilton \cite{MR592087}, a specific graph (typically an amalgamation of a complete  multigraph, possibly multipartite) is given whose edges are colored. The aim is to find conditions under which the detachment satisfies various global fairness properties (such as degree and multiplicity  conditions) as well as fairness properties within each color. For a recent survey, we refer the reader to \cite{MR3326773}. This approach is quite powerful in solving problems in two particular areas: graph decompositions, and completing partial combinatorial structures such as latin squares, edge-colored graphs, etc. \cite{MR583213,MR584118,MR1863393,MR592087,MR746544, MR920647,MR820838,MR2325799,MR1315436}.

The results of the first line of work are typically interesting from a pure mathematical point of view, while the second line of research is more appealing to those with an applied viewpoint. The goal of this paper is to study detachments by combining these two approaches. Unlike most results in the literature where the focus is graphs, here we work with  arbitrary hypergraphs whose edges are colored. 
Our goal is to find a detachment  ensuring connectivity of each color class is maintained and at the same time several fairness properties are satisfied in the resulting hypergraph as a whole as well as within each color class.  

In order to state our main result, we need a few pieces of notation. More precise definitions are provided in Section \ref{notationsec}. 

Here $x\approx y$ means $\lfloor y \rfloor \leq x\leq \lceil y \rceil$. For a hypergraph $\cG:=(V,E)$ whose edges are colored with $k$ colors (the set of colors is always $[k]:=\{1,\dots,k\}$), 
 $c(\cG)$ is the number of components of $\cG$; for $j\in [k]$, $\cG(j)$ is the color class $j$ of $\cG$; for $\alpha\in V$, $\dg_\cG (\alpha)$ is the total number of occurrences of $\alpha$ among all edges of $E$; $\omega_\al(\cG)$ is the number of ``pieces" of $\cG$ after removal of the vertex $\alpha\in V$ (our imprecise definition here will be followed by a precise one in Section \ref{wingsec}); and for $t\geq 0$, $U,X\subset V$,  $\mult_\cG(\alpha^{t},X)$ and $\mult_\cG(U,X)$ are the number of repetitions of $\{\alpha^{t}\}\cup X$ and $U\cup X$ in $E$, respectively. Finally, $\cG$ is called {\it $\alpha$-simple} if $\mult_\cG(\alpha^{2},X)=0$ for every $X\subset V$.

Here is our main result.
\begin{theorem} \label{mainthconfdet} 
Let $k$ and $n$ be positive integers. Let $\cH:=(V,E)$ be a hypergraph whose edges are colored with $k$ colors, and let $\alpha\in V$ such that 
\begin{align*}
\mult_{\cH}(\alpha^t,Y)=0 \mbox{ for each } t>n \mbox { and each }Y\subseteq V.
\end{align*}

Then there exists  an $(\alpha,n)$-detachment $\cF$ of $\cH$ such that $\cF$ is $\alpha$-simple and  the following conditions hold. 
\begin{enumerate}
\item [\textup{(F1)}] For each $i\in [n]$,
\begin{align*}
\dg_\cF (\alpha_i)&\approx \frac{\dg_\cH(\alpha)}{n};
\end{align*}
\item [\textup{(F2)}] For each $i\in [n]$ and each $j\in [k]$,
\begin{align*}
\dg_{\cF(j)} (\alpha_i)&\approx \frac{\dg_{\cH(j)}(\alpha)}{n};
\end{align*}
\item [\textup{(F3)}] For each $U\subseteq \{\alpha_1,\dots,\alpha_n\}$  and each $X\subseteq V\backslash\{\alpha_1,\dots,\alpha_n\}$, 
\begin{align*}
\mult_\cF(U,X)&\approx  \frac{\mult_\cH(\alpha^{|U|},X)}{\binom{n}{|U|}};
\end{align*}
\item [\textup{(F4)}] For each $U\subseteq \{\alpha_1,\dots,\alpha_n\}$, each $X\subseteq V\backslash\{\alpha_1,\dots,\alpha_n\}$,  and each $ j\in [k]$,
\begin{align*}
\mult_{\cF (j)}(U,X)&\approx  \frac{\mult_{\cH(j)}(\al^{|U|},X)}{\binom{n}{|U|}};
\end{align*}
\item [\textup{(C1)}] For each $j\in [k]$,
\begin{align*}
c(\cF(j))=c(\cH(j)) \mbox { if and only if } \dg_{\cH(j)}(\al)-\omega_\al(\cH(j))\geq n-1.
 \end{align*}
\end{enumerate}
\end{theorem}
We refer to conditions (F1)--(F4), and (C1) as fairness, and connectivity  condition(s), respectively. 

A special case of Theorem \ref{mainthconfdet} was known for the case when $\cH$ is an amalgamation of the  complete $n$-vertex graph \cite{MR1983358} (Hilton, Johnson,  Rodger, and Wantland \cite[p.~268]{MR1983358} noted that  the case of connectivity is ``a tough nut to crack"). For other results on connected detachments of graphs, we refer the reader to \cite{MR2946077,MR3070055,MR2325799,MR916377,MR1315436}. For detachments of 3-uniform hypergraphs and arbitrary hypergraphs, see  \cite{MR2989808} and \cite{MR2942724}, respectively (neither of these two papers address connectivity). We also note that previously, Theorem \ref{mainthconfdet} was not even known for the  case when $\cH$ is an arbitrary graph. Moreover, necessary and sufficient conditions under which a hypergraph has a connected $(\alpha,n)$-detachment were not known (condition (C1) by itself  settles this).  

The paper is organized as follows.  More precise definitions are given in Section \ref{notationsec}. In Section \ref{wingsec}, we discuss all the necessary tools needed to deal with connectivity in Theorem \ref{mainthconfdet}. In Section \ref{lamfamdefsec}, we construct two laminar families; these are crucial to deal with fairness conditions in Theorem \ref{mainthconfdet}. Before we delve into proving our main result, in Section \ref{connbaran}, we prove a generalization of the famous Baranyai's theorem \cite{MR0416986}. We prove our main result in Section \ref{mainthsecdet}. Finally, in Section \ref{corrfacgencon}, we decompose complete (not necessarily  uniform) hypergraphs into (almost) regular connected spanning sub-hypergraphs. The results of this section extend the Berge-Johnson theorem \cite{MR0498254}. In part II of this paper, we shall study  the problem of extending edge-colorings of hypergraphs.

\section{More Precise Notation}\label{notationsec}
In this paper, $\mathbb R, \mathbb{N}$, and $[k]$ denote the set of reals, positive integers,  and  $\{1,\ldots,k\}$, respectively.  For $x, y\in \mathbb{R}$, $\lfloor x \rfloor$ and $\lceil x \rceil$ denote the integers such that $x-1<\lfloor x \rfloor \leq x \leq \lceil x \rceil < x+1$, and $x\approx y$ means $\lfloor y \rfloor \leq x\leq \lceil y \rceil$.  The relation $\approx$ is transitive, and for $x, y \in \mathbb{R}, a, b, c \in \mathbb{Z},n\in \mathbb{N}$: (i) $x\approx y$ implies $x/n\approx y/n$; (ii) $a=b-c$ and $c\approx x$, implies $a\approx b-x$.
These properties of $\approx$ will be used when required without further explanation. For $x\in \mathbb{R}$, $\frc(x)$ denotes the fractional part of $x$ which is $x-\lfloor x \rfloor$. 

A {\it hypergraph} $\mathcal G$ is a pair $(V(\mathcal G),E(\mathcal G))$ where $V(\mathcal G)$ is a finite set called the {\it vertex} set, $E(\mathcal G)$ is the {\it edge} multiset, where every edge is itself a multi-subset of $V(\mathcal G)$. This means that not only can an edge  occur multiple times in $E(\mathcal G)$, but also each vertex can have multiple occurrences within an edge. A hypergraph is said to be {\it non-trivial} if it has at least one edge. By an edge of the form $\{u_1^{m_1},u_2^{m_2},\dots,u_s^{m_s}\}$ (or for short a $u_1^{m_1}\dots u_s^{m_s}$-edge), we mean an edge in which vertex $u_i$ occurs $m_i$ times for $1\leq i\leq r$, and no other vertices are present. The \textit{size} of a $u_1^{m_1}\dots u_s^{m_s}$-edge is $m_1+\dots m_s$. 
The total number of occurrences of a  vertex $v$ among all edges of $E(\mathcal G)$ is called the {\it degree}, $\dg_{\mathcal G}(v)$ of $v$ in $\mathcal G$, and $\Delta(\cG)$ denotes the maximum degree of $\cG$. 
For an edge $e$ in  $\mathcal G$, $|e|$ is the number of distinct vertices in $e$, and the {\it multiplicity} of  $e$ , written $\mult_{\mathcal G} (e)$, is the number of repetitions of $e$ in $E(\mathcal G)$ (note that $E(\mathcal G)$ is a multiset, so an edge may appear multiple times). If $\{u_1^{m_1},\dots,u_s^{m_s}\}$ is an edge in $\mathcal G$, then we abbreviate $\mult_{\mathcal G}(\{u_1^{m_1},\dots,u_s^{m_s}\})$ to $\mult_{\mathcal G}(u_1^{m_1}\dots u_s^{m_s})$. 
If $U_1,\dots,U_s$ are multi-subsets of $V(\mathcal G)$, then $\mult_{\mathcal G}(U_1,\dots,U_s)$ means $\mult_{\mathcal G}(\bigcup_{i\in[s]} U_i)$, where the union of $U_i$s is the usual union of multisets. Whenever it is not ambiguous, we drop the subscripts; for example we write $\dg(v)$ and $\mult(e)$ instead of $\dg_{\mathcal G}(v)$ and $\mult_{\mathcal G} (e)$, respectively.

We say that  a $u_1^{m_1}\dots u_s^{m_s}$-edge possesses $m_1+\dots+m_s$ objects called \textit{hinges} such that $m_i$ of these hinges are incident to vertex $u_i$ for $1\leq i\leq s$.

A {\it $k$-edge-coloring} of $\mathcal G$ is a mapping $f:E(\mathcal G)\rightarrow [k]$  and color class $i$ of $\mathcal G$, written $\mathcal G(i)$, is the spanning sub-hypergraph of $\mathcal G$ induced by the edges of color $i$.

Let $\mathcal G$ be a hypergraph, let $U$ be
some finite set, and let $\Psi : V(\mathcal G) \to U$ be a surjective mapping.  The map $\Psi$ extends
naturally to $E(\mathcal G)$.  For $A \in E(\mathcal G)$ we define $\Psi(A) = \{\Psi(x) : x \in A\}$.
Note that $\Psi$ need not be injective, and $A$ may be a multiset.  Then we define the
hypergraph $\mathcal F$ by taking $V(\mathcal F)=U$ and $E(\mathcal F)=\{ \Psi(A) : A \in E(\mathcal G) \}$.  We say
that $\mathcal F$ is an {\it amalgamation} of $\mathcal G$, and that $\mathcal G$ is a {\it detachment} of
$\mathcal F$.  Associated with $\Psi$ is a (number) function $g$ defined by $g(u)=|\Psi^{-1}(u)|$; to be
more specific we will say that $\mathcal G$ is a $g$-detachment of $\mathcal F$.  Then $\mathcal G$ has $\sum_{u\in
  V(\mathcal F)} g(u)$ vertices.  Note that $\Psi$ induces a bijection between the edges of $\mathcal F$ and
the edges of $\mathcal G$, and that this bijection preserves the size of an edge.  We adopt the
convention that it preserves the color also, so that if we amalgamate or detach an
edge-colored hypergraph the amalgamation or detachment preserves the same coloring on the
edges.  

Intuitively speaking, an {\it $(\alpha,n)$-detachment} of $\mathcal G$ is a hypergraph obtained by splitting a vertex $\alpha$ into $n$ vertices, say $\alpha_1,\dots,\alpha_n$, and sharing the incident hinges   among the subvertices. In an $\alpha$-detachment $\mathcal G'$ of $\mathcal G$ in which we split $\alpha$ into $\alpha$ and $\beta$,  an edge of the form $\{\alpha^p,u_1,\dots,u_z\}$ in $\mathcal G$ will be of the form $\{\alpha^{p-i},\beta^{i},u_1,\dots,u_z\}$ in $\mathcal G'$ for some $i$, $0\leq i\leq p$. Note that a hypergraph and its detachments have the same hinges.

If we replace every edge $e$ of a hypergraph $\cc G$ by  $\lambda$ copies of $e$, then we denote the new hypergraph by $\lambda \cc G$. For hypergraphs $\cc G_1,\dots, \cc G_t$ with the same vertex set $V$, we define their {\it union}, written $\bigcup_{i\in[t]}\cc G_i$, to be  the hypergraph with vertex set $V$ and edge set $\bigcup_{i\in[t]}E(\cc G_i)$. For a hypergraph $\cc G$ and $V\subseteq V(\cc G)$, let $\cc G-V$ be the hypergraph whose vertex set is $V(\cc G)\backslash V$ and whose edge set is $\{e\backslash V|e\in E(\cc G)\}$, and let $\cc G\backslash V$ be the hypergraph whose vertex set is $V(\cc G)$ and whose edge set is $\{e\in E(\cc G)|e\nsubseteq V\}$. 

A hypergraph $\scr G$ is said to be \textit{$k$-regular}  if every vertex has degree $k$. A \textit{$k$-factor} of a hypergraph $\scr G$ is a $k$-regular spanning sub-hypergraph of $\scr G$. 

\section{Wings}\label{wingsec}
For two hypergraphs $\mathcal G$ and $\mathcal H$, we define $\mathcal G\cup \mathcal H$ and $\mathcal G\cap \mathcal H$ as  the hypergraphs where the vertex sets are $V(\mathcal G)\cup V(\mathcal H)$ and $V(\mathcal G)\cap V(\mathcal H)$, respectively, and the edge sets are  $E(\mathcal G)\cup E(\mathcal H)$ and $E(\mathcal G)\cap E(\mathcal H)$, respectively. 
A vertex $\alpha$ in a connected hypergraph $\mathcal G$ is a {\it cut vertex} if  
there exist two non-trivial sub-hypergraphs $I, J$ of $\mathcal G$ such that 
\begin{enumerate} [label=(\roman*)]
\item $I\cup J=\mathcal G$,
\item $V( I\cap  J)=\{\alpha\}$, and 
\item $E( I\cap  J)=\varnothing$.
\end{enumerate}
A sub-hypergraph $W$ of a  hypergraph $\mathcal G$ is an {\it $\alpha$-wing} of $\mathcal G$ if 
\begin{enumerate} [label=(\roman*)]
\item $W$ is non-trivial and connected,
\item $\alpha\in V(W)$ and $\alpha$ is not a cut vertex of $W$, and 
\item no edge in $E(\mathcal G)\backslash E(W)$ is incident with a vertex in $V(W)\backslash \{\alpha\}$. 
\end{enumerate}
We remark that the edge sets of the $\alpha$-wings partition the edge set of the component of $\mathcal{G}$ containing $\alpha$.   The set of all $\alpha$-wings of $\mathcal G$ is denoted by $\scr W_\alpha (\mathcal G)$, and $\omega_{\alpha}(\cG):=|\scr W_\alpha (\mathcal G)|$. 
 If $\mathcal G$ is non-trivial and connected and $\alpha$ is not a cut vertex of $\mathcal G$, then $\scr W_\alpha (\mathcal G)=\{\mathcal G\}$. Moreover, if $\cG$ is disconnected, then $\scr W_\alpha (\mathcal G)=\scr W_\alpha (\mathcal H)$ where $\cH$ is the component of $\cG$ that contains $\alpha$. 

 An $\alpha$-wing $W$ is (i) {\it slim} if $\dg_W(\alpha)=1$, (ii) {\it wide} if $\dg_W(\alpha)\geq 2$, (iii) {\it small} if $V(W)=\{\alpha\}$, and (iv) {\it large} if $V(W)\neq \{\alpha\}$. Let $\omega^L_\alpha(\cc G)$ be the number of large $\alpha$-wings in $\cG$.

Let $\cG'$ be an $(\alpha,2)$-detachment of $\cG$, and let $W$ be an $\alpha$-wing in $\cG$. Let $W'$  denote the sub-hypergraph of $\mathcal G'$ whose hinges are the same as those in $W$. Our convention of what $W'$ refers to will continue to be used throughout this section. 
With respect to this detachment,  $W$ is (i) {\it  attached} if $\dg_{W'}(\alpha)=\dg_W(\alpha)$, (ii) {\it semi-detached} if $1\leq \dg_{W'}(\alpha)<\dg_W(\alpha)$, and (iii) {\it  detached} if $\dg_{W'}(\alpha)=0$. We remark that this terminology only applies when one of the new vertices in the detachment is called $\alpha$. Note that semi-detached wings are always wide.

\begin{lemma}\label{unionwingsimp}
Let $\cG, \cH$ be two connected hypergraphs such that $V(\cG\cap \cH)=\{\beta\}$ and $E(\cG\cap \cH)=\varnothing$. 
If $W$ is a large $\alpha$-wing of $\cG$, and $\alpha \neq\beta\in V(W)$, then $W\cup \mathcal H$ is an $\alpha$-wing of $\cG \cup \cH$. 
\end{lemma}
\begin{proof}
Since $W$ is an $\alpha$-wing, $W\cup\cH$ is non-trivial. Since both $\cG$ and $\cH$ are connected, and $V(\cG)\cap V(\cH)\neq \varnothing$, $\cG\cup\cH$ is also connected. 
Moreover, $\alpha$ is not-cut vertex of $W$ in $\mathcal G$, so it is also a non-cut vertex of $W \cup \cH$ in $\cG \cup \cH$. It is also clear that no edge in $E(\mathcal G\cup \mathcal H)\backslash E(W\cup \mathcal H)=E(\mathcal G)\backslash E(W)$ is incident with a vertex in $V(W\cup \mathcal H)\backslash \{\alpha\}$. 
Therefore, $W\cup \mathcal H$ is an $\alpha$-wing of $\cG \cup \cH$.  
\end{proof}

\begin{lemma} \label{wingplusjunk}
Let $\cG'$ be an $(\alpha,2)$-detachment of a hypergraph $\cG$. Let $X_1,\dots,X_r$ be the list of all detached $\alpha$-wings, and let $W_1,\dots,W_s$ be the  list of all semi-detached $\alpha$-wings. If $s\geq 1$, then
$$
Z:=(\cup_{i\in [r]} X'_i)\bigcup( \cup_{i\in[s]} W'_i)
$$
 is an $\alpha$-wing of $\cG'$. 
\end{lemma}
\begin{proof}
Let $X=\{X_1,\dots,X_r\}$.
 There are two cases to consider.
\begin{enumerate}
\item $X=\varnothing$: Since $s\geq 1$, $Z$ is non-trivial. For  all $i\in [s]$, $W_i$ is connected and  $\alpha\in V(W_i)$, and so $Z$ is  connected. Let $\{\beta\}=V(\cG')\backslash V(\cG)$.  Since $W_i$ is semi-detached, $\{\alpha, \beta\}\subseteq V(W'_i)$ for all $i\in [s]$.  First we show that no edge in $E(\mathcal G')\backslash E(Z)$ is incident with a vertex in $V(Z)\backslash \{\alpha\}$.  Suppose to the contrary that there is some edge $e\in E(\mathcal G')\backslash E(Z)$ with  $e\cap (V(Z)\backslash \{\alpha\})\neq \varnothing$. Since $\beta\notin e$, we have
$$
e\in E(\mathcal G)\backslash E(W_i) \mbox { and } e\cap (V(W_i)\backslash \{\alpha\})\neq \varnothing \mbox{ for some }i\in[s],
$$
 contradicting the fact that $W_i$ is an $\alpha$-wing in $\cG$. 

To complete the proof, we need to show that $\alpha$ is not a cut vertex of $Z$. Suppose to the contrary  that $\alpha$ is a cut vertex of $Z$. There must exist two non-trivial sub-hypergraphs $I', J'$ of $Z$ such that $$I'\cup J'=Z, V( I'\cap  J')=\{\alpha\}, E( I'\cap  J')=\varnothing.$$
  Without loss of generality, let us assume that $\beta \in V(I')$ (so $\beta \notin V(J')$). For each $i\in [s]$, define $I'_i=W_i'\cap I', J'_i=W_i'\cap J'$. We have  
 \begin{align*}
&I_i'\cup J_i'= (W_i'\cap I')\cup (W_i'\cap J')=W_i'\cap (I'\cup J')=W_i'\cap \cG'=W_i',\\
&V(I_i'\cap J_i')= V(I'\cap J'\cap W_i')=\{\alpha\},\\
&E(I_i'\cap J_i')\subseteq  E(I'\cap J')=\varnothing.
  \end{align*}
If for some $i\in [s]$, $I'_i$ is trivial, then $J'_i=W_i'$. But $\beta\notin V(J')$ and $\beta \in V(W_i')$, so this is clearly a contradiction.  Thus all $I_i'$s are non-trivial.  Moreover, since $J'$ is non-trivial and $J'=\cup_{i=1}^s J_i'$, there is some $k\in [s]$ so that $J'_k$ is non-trivial. Now, let $I_k$ and $J_k$ be the sub-hypergraphs of $\cG$ obtained by amalgamating $\alpha$ and $\beta$ in $I'_k$ and $J'_k$, respectively. Both $I_k$ and $J_k$ are non-trivial sub-hypergraphs of $W_k$, and 
$$I_k\cup J_k=W_k, V(I_k\cap J_k)=\{\alpha\}, E(I_i\cap J_k)=\varnothing.$$
    So, $\alpha$ is a cut vertex of $W_k$, which contradicts  $W_k$ being an $\alpha$-wing. 
\item $X\neq \varnothing$: Combining the first case with Lemma \ref{unionwingsimp} settles this case. 
\end{enumerate}
\end{proof}
\begin{lemma} \label{3partwinglem1}
Let $\cG'$ be an $(\alpha,2)$-detachment of $\cG$ with $V(\cG')=V(\cG)\cup \{\beta\}$, and let $A$ and $B$ be the set of those $\alpha$-wings of $\cG$ that are  detached and semi-detached in $\cG'$, respectively. Then
\begin{enumerate}[label=\textup{(\roman*)}]
\item $\omega_\alpha(\cG')=\omega_\alpha(\cG)-|A|-|B|+1$; 
\item $\dg_{\cG'}(\alpha)=\dg_\cG(\alpha)-\sum_{W\in A}\dg_W(\alpha)-\sum_{W\in B} \dg_{W'}(\beta)$;
\item $\dg_{\cG'}(\alpha)-\omega_\alpha(\cG')\leq \dg_\cG(\alpha)-\omega_\alpha(\cG)-1$.
\end{enumerate}
\end{lemma}
\begin{proof}
The proof of (i) is an immediate consequence of the previous lemma, and  (ii) is sufficiently obvious. For every $W\in A$, $\dg_W(\alpha)\geq 1$ and for every $W\in B$, $\dg_{W'}(\beta)\geq 1$. By combining (i) and (ii), we have
\begin{align*}
\dg_{\cG'}(\alpha)-\omega_\alpha(\cG')&=\dg_\cG(\alpha)-\sum_{W\in A}\dg_W(\alpha)-\sum_{W\in B} \dg_{W'}(\beta)-\big(\omega_\alpha(\cG)-|A|-|B|+1\big)\\
&=\dg_\cG(\alpha)-\omega_\alpha(\cG)-\sum_{W\in A}\big(\dg_W(\alpha)-1\big)-\sum_{W\in B}\big( \dg_{W'}(\beta)-1\big)-1\\
&\leq   \dg_\cG(\alpha)-\omega_\alpha(\cG)-1.
\end{align*}
This proves (iii).
\end{proof}
\begin{lemma} \label{simpleconlem}
Let  $\mathcal G'$ be an $(\alpha,2)$-detachment of a  hypergraph $\cG$ with $V(\cG')=V(\cG)\cup\{\beta\}$. Then $c(\mathcal G')=c(\cG)$  if and only if some wide $\alpha$-wing in $\cG$ is semi-detached in $\cG'$.
\end{lemma}
\begin{proof}
This is immediate from definitions. 
\end{proof}  
\begin{lemma} \label{neccconlemimp}
Let $\cF$ be an $(\alpha,n)$-detachment of a hypergraph $\cG$ such that $c(\cF)=c(\cG)$. Then
$$
\dg_\cG(\alpha)-\omega_\alpha(\cG)\geq n-1. 
$$
\end{lemma}
\begin{proof}
The proof is by induction on $n$. The case $n=1$ is trivial, because for each $W\in \scr W_\alpha (\mathcal G)$, $\dg_W(\alpha)\geq 1$. The case $n=2$ is immediate from Lemma \ref{simpleconlem}. 

To prove the inductive step, let $\cF$ be an $(\alpha,n)$-detachment of $\cG$ obtained by splitting $\alpha$ into $\alpha_1,\dots,\alpha_n$ such that $c(\cF)=c(\cG)$. Let $\cH$ be the hypergraph obtained by amalgamating $\alpha_2,\dots, \alpha_n$ in $\cF$ into $\alpha_2$ in $\cH$.  Since $\cH$ is an $(\alpha,2)$-detachment of $\cG$, by Lemma \ref{3partwinglem1}(iii) we have
$$
\dg_\cG(\alpha)-\omega_\alpha(\cG)\geq \dg_\cH(\alpha)-\omega_\alpha(\cH)+1.
$$
Since $\cF$ is an $(\alpha,n-1)$-detachment of $\cH$ and $c(\cH)=c(\cF)$, by the induction hypothesis we have
$$
\dg_\cH(\alpha)-\omega_\alpha(\cH)\geq n-2.
$$
Therefore, $\dg_\cG(\alpha)-\omega_\alpha(\cG)\geq n-1$. 
\end{proof}

\section{Laminar Families of Hinges}   \label{lamfamdefsec}
A family $\scr A$ of sets is \textit{laminar} if, for every pair $A, B$ of sets belonging to $\scr A$, either $A\subseteq B$, or $B\subseteq A$, or $A\cap B=\varnothing$.

For the rest of this section, $\cG=(V,E)$ is a $k$-edge-colored hypergraph, and  $\alpha\in V$. The set of hinges incident with $\alpha$ in $\cG$ will be denoted by $\mathbb H_\cG(\alpha)$. In this section, we construct two laminar families of subsets of $\mathbb H_\cG(\alpha)$. These families will be crucial in the proof of our main result. For an edge $e$ in $\cG$,  $\mathbb H_{\cG}^e(\alpha)$  will be the set of hinges incident with both $\alpha$ and $e$ in $\cG$; abbreviate $\mathbb H_\cG^{\{\alpha^r\}\cup U}(\alpha)$ to $\mathbb{H}_\cG^{r,U}(\alpha)$ (for $r\geq 1$ and $U \subseteq V\backslash\{\alpha\}$). 
 For $i\in [k]$, let 
 \begin{align*}
 ^i \mathbb H_\cG(\alpha)=\bigcup_{W\in \scr W_\alpha(\mathcal G(i)), \dg_W(\alpha)\geq 2} \mathbb H_W(\alpha).
\end{align*}
In other words, $^i \mathbb H_\cG(\alpha)$ is the set of all hinges incident with $\alpha$ in some wide $\alpha$-wing of $\cG(i)$.
Observe that 
\begin{itemize}
\item $|\mathbb{H}_\cG(\alpha)|=\dg_\cG(\alpha)$,
\item $|\mathbb{H}_{\cG(i)}(\alpha)|=\dg_{\cG(i)}(\alpha)$ for  $i\in[k]$,
\item $\mathbb H_{\cG}^e(\alpha)\subseteq \mathbb H_{\cG(i)}(\alpha)$ for $i\in[k]$ and $e\in E(\cG(i))$,
\item $|\mathbb H_\cG^{r,U}(\alpha)|=r\mult_\cG(\alpha^r,U)$ for $r\geq 1$ and $U \subseteq V\backslash\{\alpha\}$. 
\item $|\mathbb H_{\cG(i)}^{r,U}(\alpha)|=r\mult_{\cG(i)}(\alpha^r,U)$ for $i\in [k], r\geq 1$ and $U \subseteq V\backslash\{\alpha\}$. 
\end{itemize}
 
 The proof of the following lemma is immediate from definitions. 
\begin{lemma}
Let
 \begin{align*}
\scr A  &= \quad \mathbb H_\cG(\alpha)\\
&\quad  \bigcup    \{\mathbb H_{\cG(1)}(\alpha),\ldots,\mathbb H_{\cG(k)}(\alpha)\} \\
&\quad \bigcup    \{^1 \mathbb H_\cG(\alpha),\ldots,{^k \mathbb H_\cG(\alpha)}\} \nonumber\\
&\quad \bigcup  \{ \mathbb H^e_\cG(\alpha) : e  \in E\}\\
&\quad  \bigcup  \{ \mathbb H_{W}(\alpha) : W  \in \scr W_{ \alpha}(\mathcal G(i)), \dg_W(\alpha)\geq 2, i \in [k] \}   \nonumber,
\end{align*}
and 
\begin{align*} 
\scr B\quad &= \{\mathbb H_\cG^{r, U}(\alpha): r\geq 1, U \subseteq V\backslash\{\alpha\}\}\\
&\quad \bigcup \{\mathbb H^{r, U}_{\cG(i)}(\alpha): r\geq 1, U \subseteq V\backslash\{\alpha\}, i \in [k]\}.
\end{align*}
Then both $\scr A$ and $\scr B$  are  laminar families of subsets of $\mathbb H_\cG(\alpha)$. 
\end{lemma} 

\section{Connected Baranyai's Theorem} \label{connbaran}
A hypergraph is {\it almost regular} if there is an integer $d$ such that the degree of each vertex is $d$ or $d+1$. We exhibit the usefulness of our main result by providing the following applications. 
\begin{theorem} \label{connbarayalm}
Let $a_1,\dots,a_k\in \mathbb{N}$ and $\sum_{j\in [k]} a_j=\lambda\binom{n}{h}$. Then the complete $\lambda$-fold $n$-vertex $h$-uniform hypergraph $\lambda K_n^h$ can be decomposed into almost regular spanning sub-hypergraphs $\cF(1), \dots, \cF(k)$ so that for all $j\in[k]$ the following hold.
\begin{enumerate}[label=\textup{({\roman*})}] 
\item $\cF(j)$ has $a_j$ edges;
\item $\dg_{\cF(j)} (v)\approx ha_j/n$ for each $v\in V(\cG_j)$; 
\item  $\cF(j)$ is connected if and only if  $a_j \geq (n-1)/(h-1)$.  
\end{enumerate}
\end{theorem}
\begin{proof}
Given $a_1,\dots,a_k$ with $\sum_{j\in [k]} a_j=\lambda\binom{n}{h}$, let $\cH$ be a hypergraph with 
\begin{align*} 
    &&
    V(\cH)=\{\alpha\},
    &&
    E(\cH)=\{\{\alpha^h\}^{\lambda\binom{n}{h}}\}.
    &&
  \end{align*}
We color the edges of $\cH$ with $k$ colors so that 
$$\mult_{\cH(j)}(\alpha^h)=a_j \mbox { for }j\in [k].$$ 
It is clear that for every $j\in [k]$,  $\dg_{\cH(j)}(\alpha)=ha_j$, $c(\cH(j))=1$, and $\omega_\al(\cH(j))=\mult_{\cH(j)}(\alpha^h)=a_j$. By Theorem \ref{mainthconfdet}, there exists  an $(\alpha,n)$-detachment $\cF$ of $\cH$ such that the following conditions hold.
\begin{enumerate}[label=(\alph*)]
\item By \textup{(F2)}, for each $i\in [n]$ and each $j\in [k]$,
\begin{align*}
\dg_{\cF(j)} (\alpha_i)&\approx \frac{ha_j}{n};
\end{align*}
\item By \textup{(F3)}, for each $U\subseteq \{\alpha_1,\dots,\alpha_n\}$ with $|U|=h$, 
\begin{align*}
\mult_\cF(U)&\approx  \frac{\mult_\cH(\alpha^{h})}{\binom{n}{h}}=\frac{\lambda \binom{n}{h}}{\binom{n}{h}}=\lambda;
\end{align*}
\item By \textup{(C1)}, for each $j\in [k]$,
\begin{align*}
\cF(j)  \mbox { is connected if and only if } ha_j-a_j\geq n-1.
 \end{align*}
\end{enumerate}
\end{proof}
Theorem \ref{connbarayalm} without condition (iii) is known as the Baranyai's theorem \cite{MR0416986}. A special case of Theorem \ref{connbarayalm} where almost regularity is replaced by regularity is proven in \cite{MR3213845}. 

\section{Proof of Theorem \ref{mainthconfdet}} \label{mainthsecdet}
Let $k,n\in \mathbb{N}$. Let $\cH=(V,E)$ be a $k$-edge-colored hypergraph with $\alpha\in V$ such that 
\begin{align*}
\mult_{\cH}(\al^t,Y)=0 \mbox{ for each } t>n \mbox { and each }Y\subseteq V.
\end{align*}
We show that there exists  an $(\alpha,n)$-detachment $\cF$ of $\cH$ such that $\cF$ is $\alpha$-simple, and that conditions (F1)--(F4) and (C1) hold.

Let $I\subseteq [k]$ be the set of all colors so that for each $i\in I$, $\dg_{\cH(i)}(\al)-\omega_\al(\cH(i))\geq n-1$. 

Given an arbitrary $(\alpha,n)$-detachment $\cF$ of $\cH$, if for some $i\in [k]$,  $c(\cH(i))=c(\cF(i))$, then by  Lemma \ref{neccconlemimp} we must have $i\in I$. Since an $\alpha$-detachment only affects the component containing $\alpha$, to prove (C1) it  suffices to show that if $\cH(i)$ is connected, then so is $\cF(i)$,  for all $i\in I$. For the rest of this section we will assume that $\cH(i)$ is connected for all $i\in I$. 

We shall construct  a sequence of hypergraphs $\cH_1,\dots,\cH_n$ with $\cH_1:=\cH$, $\cH_n:=\cF$, $\alpha_1:=\alpha$, and $V(\cH_\ell)=V(\cH_{\ell-1})\cup\{\alpha_{\ell}\}$ for $2\leq \ell\leq n$, where $\cH_\ell$ is a detachment of $\cH_{\ell-1}$ obtained by splitting off the new vertex $\alpha_{\ell}$ from $\al$ so that we end up with $\alpha$ and $\alpha_{\ell}$ in $\cH_{\ell}$,  such that for all $\ell\in [n]$, $\cH_\ell$ satisfies the following conditions.

For all $i\in [\ell], j\in [k], 0\leq r\leq n-\ell+1, t>n-\ell+1, s\in I, Y\subseteq V, U\subseteq \{\alpha_2,\dots,\alpha_\ell\}, X\subseteq V\backslash\{\alpha_1,\dots,\alpha_\ell\}, 1\leq r+|U|\leq n$:
\begin{align} 
\mult_{\cH_\ell}(\alpha^t,Y)&=0; \label{simpleco}\\
\dg_{\cH_\ell}(\alpha_i) &\approx \left \{ \begin{array}{ll} 
\dfrac{(n-\ell+1)\dg_\cH(\alpha)}{n} & \mbox { for } i = 1  \\
\dfrac{\dg_\cH(\alpha)}{n} & \mbox { for } 2\leq i\leq \ell \end{array} \right.\label{degco};\\
\dg_{\cH_\ell(j)}(\alpha_i) &\approx \left \{ \begin{array}{ll}
\dfrac{(n-\ell+1)\dg_{\cH(j)}(\alpha)}{n} & \mbox { for } i = 1  \\
\dfrac{\dg_{\cH(j)}(\alpha)}{n} & \mbox { for } 2\leq i\leq \ell \end{array} \right. \label{degjco};\\
\frac{\mult_{\cH_\ell}(\alpha^r,U,X)}{\binom{n-\ell+1}{r}}&\approx  \frac{\mult_\cH(\alpha^{r+|U|},X)}{\binom{n}{r+|U| }};\label{multco}\\
\frac{\mult_{\cH_\ell(j)}(\alpha^r,U,X)}{\binom{n-\ell+1}{r}}&\approx  \frac{\mult_{\cH(j)}(\alpha^{r+|U|},X)}{\binom{n}{r+|U| }}\label{multjco};\\
\cH_\ell(s) \mbox{ is connected;} \label{connco}\\
\dg_{\cH_\ell(s)}(\al)-\omega_\al(\cH_\ell(s))&\geq n-\ell \label{necineqco}.
\end{align}

We prove this by induction on $\ell$. For $\ell=1$,  $\cH$ clearly satisfies all the conditions \eqref{simpleco}--\eqref{necineqco}, so we can let $\cH_1$  be $\cH$. 

Before we proceed with the inductive step, we verify that when $\ell=n$ (and so $n-\ell+1=1$), we can let $\cF=\cH_n$: Condition \eqref{simpleco} implies that $\cH_n$ is $\alpha$-simple, and conditions \eqref{degco} and \eqref{degjco} are  equivalent to (F1) and (F2) respectively. To see the equivalence of \eqref{multco} and (F3), notice that the set $U$ in (F3) either contains $\alpha_1=\alpha$ or it does not. 
According to condition \eqref{multco}, for $0\leq r\leq 1, U\subseteq \{\alpha_2,\dots,\alpha_n\}, X\subseteq V\backslash\{\alpha_1,\dots,\alpha_n\}, 1\leq r+|U|\leq n$ we have
$$
\frac{\mult_{\cH_n}(\alpha^r,U,X)}{\binom{1}{r}}\approx  \frac{\mult_\cH(\alpha^{r+|U|},X)}{\binom{n}{r+|U| }}.
$$
For $r=0$, this is equivalent to
$$
\mult_{\cH_n}(U,X)\approx  \frac{\mult_\cH(\alpha^{|U|},X)}{\binom{n}{|U| }},
$$
and for $r=1$, it is equivalent to
$$
\mult_{\cH_n}(\alpha,U,X)\approx  \frac{\mult_\cH(\alpha^{|U|+1},X)}{\binom{n}{|U|+1 }}.
$$
The proof of the equivalence of conditions \eqref{multjco} and (F4) is very similar. Condition \eqref{connco} together with the discussion at the beginning of this section implies (C1). 

To prove the inductive step, let us assume that $\cG:=\cH_\ell$ has been constructed for some $\ell<n$, and $\cH_\ell$ satisfies \eqref{simpleco}--\eqref{necineqco}. We denote the corresponding conditions for $\cG':=\cH_{\ell+1}$ by \eqref{simpleco}$'$--\eqref{necineqco}$'$.

In order to decide how to share the edges incident with $\alpha$ in $\cG$ between $\alpha$ and the new vertex $\alpha_{\ell+1}$ in $\cG'$, let  $\scr A$ and $\scr B$ be  the  laminar families of subsets of $\mathbb H_\cG(\alpha)$ defined in Section \ref{lamfamdefsec}. By \cite[Lemma 2]{MR916377}, there exists a subset $Z$ of $\mathbb{H}_\cG(\alpha)$ such that 
\begin{align} \label{lamnashcon}
 |Z\cap P|\approx \frac{|P|}{n-\ell+1} \mbox { for every } P\in \scr A \cup \scr B. 
\end{align} 
Let $\mathcal G'$ be the $(\alpha,2)$-detachment of $\mathcal G$ with $V(\cG')\backslash V(\cG)=\{\beta\}:=\{\alpha_{\ell+1}\}$
such  that the hinges which were incident with $\alpha$ in $ \mathcal G$ are incident in $\mathcal G'$ with $\alpha$ or $\beta$ according to whether  they do not or do belong to $Z$, respectively. More precisely,
 \begin{equation}\label{hinge1'} 
 \mathbb H_{\cG'}(\beta)=Z, \quad \mathbb H_{\cG'}(\alpha)=\mathbb H_\cG( \alpha)\backslash Z. 
 \end{equation}
Note that \eqref{hinge1'} provides  a full description of $\cG'$. It remains to show that $\cG'$ indeed satisfies \eqref{simpleco}$'$--\eqref{necineqco}$'$. Throughout the rest of our proof, we shall repeatedly use properties of $\approx$ and  \eqref{simpleco}--\eqref{hinge1'} without further explanation. 

By the induction hypothesis, for $t>n-\ell+1$ and  $Y\subseteq V$ we have  $\mult_\cG(\alpha^t,Y)=0$, or equivalently,   for  $e\in E(\cG)$ with $\alpha\in e$, we have $|\mathbb{H}_\cG^e(\alpha)|\leq n-\ell+1$. Since $\mathbb H^e_\cG(\alpha)\in \mathscr{A}$, we have 
$$|Z\cap \mathbb{H}^e_\cG(\alpha)|\approx \frac{|\mathbb{H}_\cG^e(\alpha)|}{n-\ell+1}\leq 1.
$$
Therefore,  $|Z\cap \mathbb{H}^e_\cG(\alpha)|\in \{0,1\}$, and so $\mult_{\cG'}(\beta^q,U)=0$ for $q\geq 2$ and  $U\subseteq V(\cG')$. 
If $|\mathbb{H}_\cG^e(\alpha)|=n-\ell+1$, then $|Z\cap \mathbb{H}_\cG^e(\alpha)|=1$, so $|\mathbb{H}_{\cG'}^e(\alpha)|=n-\ell$, and if $|\mathbb{H}_\cG^e(\alpha)|<n-\ell+1$, then $|\mathbb{H}_{\cG'}^e(\alpha)|\leq |\mathbb{H}^e_\cG(\alpha)|\leq n-\ell$. Therefore, in both cases \eqref{simpleco}$'$ is satisfied.

Since $\mathbb H_\cG(\alpha)\in \scr A$, we have  
\begin{align*}
\dg_{\cG'}(\beta)& =  |Z|=|Z\cap \mathbb H_\cG(\alpha)|\\
& \approx  \frac{|\mathbb H_\cG(\alpha)|}{n-\ell+1}=\frac{\dg_\cG(\alpha)}{n-\ell+1} \\
& \approx \frac{\dg_\cH(\alpha)}{n},\\
\dg_{\cG'}(\alpha) & =   |\mathbb{H}_{\cG'}(\alpha)|=|\mathbb H_\cG( \alpha)\backslash Z|\\
&= \dg_\cG(\alpha) - \dg_{\cG'}(\beta) \nonumber \\
& \approx   \frac{(n-\ell+1)\dg_\cH(\alpha)}{n}-\frac{\dg_\cH(\alpha)}{n}\\
&= \frac{(n-\ell)\dg_\cH(\alpha)}{n}.
\end{align*}
Moreover, for $2\leq i\leq \ell$ we have $\dg_{\cG'}(\alpha_i)=\dg_\cG(\alpha_i)=\dg_{\cH}(\alpha)/n$. 
Therefore $\mathcal G'$ satisfies \eqref{degco}$'$. 

For $j\in [k]$, we have $\mathbb H_{\cG(j)}(\alpha)\in \scr A$, and so,   
\begin{align*}
\dg_{\cG'(j)}(\beta)& =  |Z\cap \mathbb H_{\cG(j)}(\alpha)|\\
& \approx  \frac{|\mathbb H_{\cG(j)}(\alpha)|}{n-\ell+1}=\frac{\dg_{\cG(j)}(\alpha)}{n-\ell+1} \\
& \approx \frac{\dg_{\cH(j)}(\alpha)}{n},\\
\dg_{\cG'(j)}(\alpha) & =    \dg_{\cG(j)}(\alpha) -\dg_{\cG'(j)}(\beta) \nonumber \\
& \approx  \frac{(n-\ell+1)\dg_{\cH(j)}(\alpha)}{n}-\frac{\dg_{\cH(j)}(\alpha)}{n}\\
&= \frac{(n-\ell)\dg_{\cH(j)}(\alpha)}{n}.
\end{align*}
This proves \eqref{degjco}$'$.

For $r\geq 1, U\subseteq V(\cG)\backslash \{\alpha\}$, we have $\mathbb{H}_\cG^{r,U}(\alpha)\in \mathscr B$. Therefore, 
\begin{align*}
\mult_{\cG'}(\alpha^{r-1},\beta,U)& =  |Z\cap \mathbb{H}_\cG^{r,U}(\alpha)|\\
& \approx  \frac{|\mathbb{H}_\cG^{r,U}(\alpha)|}{n-\ell+1}\\
&=\frac{r}{n-\ell+1}\mult_\cG(\alpha^r,U), \\
\mult_{\cG'}(\alpha^r,U) & =    \mult_\cG(\alpha^r,U) - \mult_{\cG'}(\alpha^{r-1},\beta,U)  \\
& \approx   \mult_\cG(\alpha^r,U) - \frac{r}{n-\ell+1}\mult_\cG(\alpha^r,U)\\
&= \frac{(n-\ell-r+1)}{n-\ell+1}\mult_\cG(\alpha^r,U).
\end{align*}

To prove \eqref{multco}$'$, let $U\subseteq \{\alpha_2,\dots,\alpha_{\ell+1}\}$, and let $X\subseteq V(\cG')\backslash \{\alpha_1,\dots,\alpha_{\ell+1}\}$. There are two cases to consider. 
\begin{enumerate}[label=\textup{(\roman*)}]
\item If $\beta\in U$, then
\begin{align*}
\frac{\mult_{\cG'}(\alpha^{r},U,X)}{\binom{n-\ell}{r}}& \approx  \frac{(r+1)\mult_\cG(\alpha^{r+1}, U\backslash\{\beta\}, X)}{(n-\ell+1)\binom{n-\ell}{r}}\\
& =  \frac{\mult_\cG(\alpha^{r+1}, U\backslash\{\beta\}, X)}{\binom{n-\ell+1}{r+1}}\\
&\approx  \frac{\mult_\cH(\alpha^{r+1+|U\backslash\{\beta\}|},X)}{\binom{n}{r+1+|U\backslash\{\beta\}|}}\\
&=\frac{\mult_\cH(\alpha^{r+|U|},X)}{\binom{n}{r+|U|}}.
\end{align*}
\item If $\beta\notin U$, then
\begin{align*}
\frac{\mult_{\cG'}(\alpha^{r},U,X)}{\binom{n-\ell}{r}}& \approx  \frac{(n-\ell-r+1)\mult_\cG(\alpha^{r}, U, X)}{(n-\ell+1)\binom{n-\ell}{r}}\\
& =  \frac{\mult_\cG(\alpha^{r}, U, X)}{\binom{n-\ell+1}{r}}\\
&\approx \frac{\mult_\cH(\alpha^{r+|U|},X)}{\binom{n}{r+|U|}}.
\end{align*}
\end{enumerate}
This proves \eqref{multco}$'$.

For $j\in [k], r\geq 1, U\subseteq V(\cG)\backslash \{\alpha\}$, we have $\mathbb{H}_{\cG(j)}^{r,U}(\alpha)\in \mathscr B$. Therefore, 
\begin{align*}
\mult_{\cG'(j)}(\alpha^{r-1},\beta,U)& =  |Z\cap \mathbb H^{r,U}_{\cG(j)}(\alpha)|\\
& \approx  \frac{|\mathbb H^{r,U}_{\cG(j)}(\alpha)|}{n-\ell+1}\\
&=\frac{r}{n-\ell+1}\mult_{\cG(j)}(\alpha^r,U), \\
\mult_{\cG'(j)}(\alpha^r,U) & =    \mult_{\cG(j)}(\alpha^r,U) -\mult_{\cG'(j)}(\alpha^{r-1},\beta,U)  \\
& =  \mult_{\cG(j)}(\alpha^r,U) - \frac{r}{n-\ell+1}\mult_{\cG(j)}(\alpha^r,U)\\
&= \frac{(n-\ell-r+1)}{n-\ell+1}\mult_{\cG(j)}(\alpha^r,U).
\end{align*}

To prove \eqref{multjco}$'$, let $j\in [k], U\subseteq \{\alpha_2,\dots,\alpha_{\ell+1}\}$, and let $X\subseteq V(\cG')\backslash \{\alpha_1,\dots,\alpha_{\ell+1}\}$. There are two cases to consider. 
\begin{enumerate} [label=\textup{(\roman*)}]
\item If $\beta\in U$, then
\begin{align*}
\frac{\mult_{\cG'(j)}(\alpha^{r},U,X)}{\binom{n-\ell}{r}}& \approx  \frac{(r+1)\mult_{\cG(j)}(\alpha^{r+1}, U\backslash\{\beta\}, X)}{(n-\ell+1)\binom{n-\ell}{r}}\\
& =  \frac{\mult_{\cG(j)}(\alpha^{r+1}, U\backslash\{\beta\}, X)}{\binom{n-\ell+1}{r+1}}\\
&\approx  \frac{\mult_{\cH(j)}(\alpha^{r+1+|U\backslash\{\beta\}|},X)}{\binom{n}{r+1+|U\backslash\{\beta\}|}}\\
&=\frac{\mult_{\cH(j)}(\alpha^{r+|U|},X)}{\binom{n}{r+|U|}}.
\end{align*}
\item If $\beta\notin U$, then
\begin{align*}
\frac{\mult_{\cG'(j)}(\alpha^{r},U,X)}{\binom{n-\ell}{r}}& \approx  \frac{(n-\ell-r+1)\mult_{\cG(j)}(\alpha^{r}, U, X)}{(n-\ell+1)\binom{n-\ell}{r}}\\
& \approx  \frac{\mult_{\cG(j)}(\alpha^{r}, U, X)}{\binom{n-\ell+1}{r}}\\
&=\frac{\mult_{\cH(j)}(\alpha^{r+|U|},X)}{\binom{n}{r+|U|}}.
\end{align*}
\end{enumerate}
This proves \eqref{multjco}$'$.

To prove the connectivity conditions \eqref{connco}$'$ and \eqref{necineqco}$'$, let us fix a color $i\in I$ (If $I=\varnothing$, there is nothing to prove). By Lemma \ref{simpleconlem}, in order to establish  that $\cG'$ is connected, it is enough  to show that
\begin{enumerate} [label=(\alph*)]
\item There exists an $\alpha$-wing $W$ in $\cG(i)$ such that $\dg_W(\alpha)\geq 2$;
\item For every $\alpha$-wing $W$ in $\cG(i)$ with $\dg_W(\alpha)\geq 2$, we have 
$
|Z\cap \mathbb{H}_W(\alpha)|<\dg_W(\alpha);
$
\item 
$
|Z\cap {^i \mathbb H_\cG(\alpha)}|\geq 1.
$
\end{enumerate}
To prove (a), suppose to the contrary  that for every $W\in \mathscr W_\alpha(\cG(i))$, $\dg_W(\alpha)= 1$. Then,
$$\dg_{\cG(i)}(\al)-\omega_\al(\cG(i))=\sum_{W\in  \mathscr W_\alpha(\cG(i))} \big(\dg_W(\alpha) -1\big)=0.$$
But this contradicts  \eqref{necineqco}, which says $\dg_{\cG(i)}(\al)-\omega_\al(\cG(i))\geq n-\ell\geq 1$. 

To prove (b), let $W\in \mathscr W_\alpha(\cG(i))$ with $\dg_W(\alpha)\geq 2$. Since $\mathbb H_W(\alpha)\in \mathscr A$, we have
$$
|Z\cap  \mathbb H_W(\alpha)|\approx\frac{| \mathbb H_W(\alpha)|}{n-\ell+1}=\frac{\dg_W(\alpha)}{n-\ell+1}.
$$
Since $\dg_W(\alpha)\geq 2$ and $n-\ell+1\geq 2$, we have
$$
|Z\cap  \mathbb H_W(\alpha)|< \dg_W(\alpha).
$$

To prove (c), let $A_i=\{W\in \mathscr W_\alpha(\cG(i)): \dg_W(\alpha)\geq 2 \}$. By (a), $|A_i|\geq 1$. We have
\begin{align*}
n-\ell&\leq  \dg_{\cG(i)}(\al)-\omega_\al(\cG(i))\\
&=\sum_{W\in  \mathscr W_\alpha(\cG(i))} \big(\dg_W(\alpha) -1\big)\\
&=\sum_{W\in  A_i} \dg_W(\alpha) -|A_i|\\
&\leq  |{^i \mathbb H_\cG(\alpha)}|-1.
\end{align*}
Since ${^i \mathbb H_\cG(\alpha)}\in \mathscr A$,
$$
|Z\cap  {^i \mathbb H_\cG(\alpha)}|\approx \frac{|{^i \mathbb H_\cG(\alpha)}|}{n-\ell+1}\geq 1.
$$
This completes the proof of \eqref{connco}$'$, and we conclude that $\cG'(i)$ is connected for all $i\in I$.

Now, we need to show that 
$$
\dg_{\cG'(i)}(\al)-\omega_\al(\cG'(i))\geq n-\ell-1.
$$
For $W\in \mathscr W_\alpha(\cG(i))$ with $\dg_W(\alpha)\geq 2$, we have $\mathbb H_W(\alpha)\in \scr A$, and so  
\begin{align*}
\dg_{W'}(\beta)& =  |Z\cap \mathbb H_W(\alpha)|\\
& \approx  \frac{|\mathbb H_W(\alpha)|}{n-\ell+1}=\frac{\dg_W(\alpha)}{n-\ell+1}, \\
\dg_{W'}(\alpha) & =   \dg_W(\alpha) - \dg_{W'}(\beta) \nonumber \\
& \approx  \dg_W(\alpha) - \frac{\dg_W(\alpha)}{n-\ell+1}\\
&= \frac{(n-\ell)\dg_W(\alpha)}{n-\ell+1}.
\end{align*}
Let $S_i$ be the set of detached slim $\alpha$-wings of $\cG(i)$, and let $L_i$ be the set of semi-detached (wide) $\alpha$-wings of $\cG(i)$. Note that  $L_i\neq\varnothing$. We have
$$\dg_{\cG'(i)}(\alpha)=\dg_{\cG(i)}(\alpha)-|S_i|-\sum_{W\in L_i}\dg_{W'}(\beta).$$
Moreover, by Lemma \ref{3partwinglem1}(i),  
$$\omega_\alpha(\cG'(i))=\omega_\alpha(\cG( i))-|S_i|-|L_i|+1.$$
Therefore, 
\begin{align*}
\dg_{\cG'(i)}(\al)-\omega_\al(\cG'(i))+1&=\dg_{\cG(i)}(\alpha)-|S_i|-\sum_{W\in L_i}\dg_{W'}(\beta)\\
&\quad -\big(\omega_\al(\cG(i))-|S_i|-|L_i|+1\big)+1\\
&= \big(\dg_{\cG(i)}(\alpha)-\omega_\al(\cG(i))\big)-\big(\sum_{W\in L_i}\dg_{W'}(\beta)-|L_i|\big)\\
&= \big(\dg_{\cG(i)}(\alpha)-\omega_\al(\cG(i))\big)-\sum_{W\in L_i}\big(\dg_{W'}(\beta)-1\big).
\end{align*}
Let $\overline L_i=\{W\in L_i:\dg_{W'}(\beta)\geq 2\}$. There are two cases to consider.
\begin{enumerate} [label=\textup{(\Roman*)}]
\item $\overline L_i=\varnothing$: In this case, for each $W\in L_i$, $\dg_{W'}(\beta)=1$, and so $\sum_{W\in L_i}\big(\dg_{W'}(\beta)-1\big)=0$. Thus,  
$$\dg_{\cG'(i)}(\al)-\omega_\al(\cG'(i))+1=\dg_{\cG(i)}(\alpha)-\omega_\al(\cG(i))\geq n-\ell,$$ and the proof is complete.
\item $\overline L_i\neq \varnothing$: In this case, there is some $W\in L_i$ with $\dg_{W'}(\beta)\geq 2$.  Since $W$ is a semi-detached wide $\alpha$-wing, we have $\dg_W(\alpha)\geq 3$. In fact, since $\dg_{W'}(\beta)\geq 2$ and
$$\dg_{W'}(\beta)\approx \frac{\dg_W(\alpha)}{n-\ell+1},$$ 
we have  $\dg_W(\alpha)\geq n-\ell+2$. The following completes the proof.
\begin{align*}
\dg_{\cG'(i)}(\al)-\omega_\al(\cG'(i))+1 &=  \sum_{W\in \mathscr W_\alpha(\cG(i))}\big(\dg_{W}(\alpha)-1\big)-\sum_{W\in \overline L_i}\big(\dg_{W'}(\beta)-1\big)\\
&\geq  \sum_{W\in \overline L_i}\big(\dg_{W}(\alpha)-1\big)-\sum_{W\in \overline L_i}\big(\dg_{W'}(\beta)-1\big)\\
&=\sum_{W\in \overline L_i}\big(\dg_{W}(\alpha)-\dg_{W'}(\beta)\big)\\
&=\sum_{W\in \overline L_i}\dg_{W'}(\alpha)\\
&\geq \sum_{W\in \overline L_i} \dg_W(\alpha)\frac{n-\ell}{n-\ell+1}-\frac{n-\ell}{n-\ell+1}\\
&\geq    \sum_{W\in \overline L_i}\frac{(n-\ell)(\dg_{W}(\alpha)-1)}{n-\ell+1}\\
&=    \frac{n-\ell}{n-\ell+1}\sum_{W\in \overline L_i}(\dg_{W}(\alpha)-1))\\
&\geq  \frac{n-\ell}{n-\ell+1} |\overline L_i|(n-\ell+1)\\
&= |\overline L_i|(n-\ell)\\
&\geq  n-\ell.
\end{align*}
\end{enumerate}
\qed

\section{Connected Factorizations of (Non-uniform) Complete Hypergraphs} \label{corrfacgencon}  

In this section, we partition the edge set of  complete (not necessarily  uniform) hypergraphs into (almost) regular connected spanning sub-hypergraphs. First, we prove a simple observation.
\begin{lemma}  \label{elemcomplem}
Let  $\cG=(V,E)$ be a hypergraph. We have the following. 
\begin{enumerate} [label=\textup{(\alph*)}]
\item Adding  an edge $e$ to $\cG$ decreases the number of components by some $t\in \{0,|e|-1\}$. 
\item  Deleting an edge $e$ of $\cG$ increases the number of components  by some $t\in \{0,|e|-1\}$.  
\item $$c(\cG)\geq |V|-\sum_{e\in E}(|e|-1).$$ 
\item If $\cG$ is $h$-uniform, then $c(\cG)\geq |V|-(h-1)|E|$. 
\item If $\cG$ is $h$-uniform and connected, then $|E|\geq \frac{|V|-1}{h-1}$.
\end{enumerate}
\end{lemma}
\begin{proof}
Parts (a) and (b) are sufficiently obvious. Adding one-by-one the edges of $\cG$ to a trivial hypergraph with $|V|$ vertices and applying (a) each time, completes the proof of (c). Part (d) is a special case of (c), and (e) is obtained by letting $c(\cG)=1$ in (d). 
\end{proof}

For positive integer column vectors $\Lambda:=[\lambda_1\dots\lambda_m]^T$ and $H:=[h_1\dots h_m]^T$ with $1\leq h_1<h_2<\dots<h_m\leq n$, let $\Lambda K_n^H$ denote an $n$-vertex hypergraph such that there are $\lambda_i$ edges of size $h_i$ incident with every $h_i$ vertices  for $i\in [m]$. In other words,
$$\Lambda K_n^H=\bigcup\nolimits_{i\in[m]}\lambda_i K_{n}^{h_i}.$$
For a hypergraph $\cG$ and a positive integer column vector $R:=[r_1\dots r_k]^T$, an \textit{$R$-factorization} of $\cG$ is a partition (decomposition) $\{F_1,\ldots, F_k\}$ of $E(\cG)$ in which $F_i$ is an $r_i$-factor for $i\in[k]$. 
\begin{theorem} \label{r1rkfacknhlamcon}
$\Lambda K_n^H$ is $R$-factorable if and only if
\begin{equation} \label{longfatnec1}
\sum\nolimits_{i\in[m]}  \lambda_i\binom{n-1}{h_i-1}=\sum\nolimits_{i\in [k]}r_i,
\end{equation}
and there exists a non-negative integer matrix $A=[a_{ij}]_{k \times m}$ such that 
\begin{align} \label{longfatnec2}
    &&
    AH=nR,
    &&
    \sum\nolimits_{i\in[k]}a_{ij}=\lambda_j\binom{n}{h_j} \mbox{ for } j\in [m].
    &&
  \end{align}
Moreover, for each $i\in [k]$, an $r_i$-factor that contains exactly $a_{ij}$ edges of size $h_j$ for each $j\in [s]$, is connected if and only if
\begin{equation} \label{longfatnec3}
\sum\nolimits_{j\in[m]}  a_{ij}\leq n(r_i-1)+1.
\end{equation}
\end{theorem}
\begin{proof}
To prove the necessity, suppose that $\Lambda K_n^H$ is $R$-factorable where $F_i$s are the corresponding $r_i$-factors for $i\in [k]$. Since $F_i$ is $r_i$-regular and $\Lambda K_n^H$ is $\sum_{i\in [m]} \lambda_i \binom{n-1}{h_i-1}$-regular, \eqref{longfatnec1} must hold. For $i\in[k],j\in [m]$, let $a_{ij}$ be  the number of edges  of size $h_j$ (counting multiplicities) in $F_i$. Since $\{F_1,\dots,F_k\}$ is a partition of the edges of $\Lambda K_n^H$, we have
$$\sum\nolimits_{i\in[k]} a_{ij} =\lambda_j\binom{n}{h_j} \quad  \mbox{ for }j\in [m].$$
 Double counting the degree sum in $F_i$, we obtain
$$\sum\nolimits_{j\in[m]} a_{ij} h_j=nr_i\quad  \mbox{ for } i\in[k].$$
Let us fix an $i\in [k]$. If $F_i$ is connected, then by Lemma \ref{elemcomplem} we have  
$$1\geq n-\sum\nolimits_{e\in F_i}(|e|-1)=n-\sum\nolimits_{j\in[m]} a_{ij} (h_j-1)=n-r_in+\sum\nolimits_{j\in[m]} a_{ij}.$$
Therefore, \eqref{longfatnec3} holds. 

To prove the sufficiency, let $\cH$ be a hypergraph with a single vertex $\alpha$ such that 
 $$\mult_{\cH}(\alpha^{h_j})=\lambda_j \binom{n}{h_j}\quad  \mbox{ for }j\in [m].$$
Since for $j\in [m]$, $\sum\nolimits_{i\in[k]} a_{ij} =\lambda_j\binom{n}{h_j}$, we can color the edges of $\cH$ so that 
$$\mult_{\cH(i)}(\alpha^{h_j})=a_{ij} \mbox{ for }i\in[k],j\in [m].$$
Let us fix an $i\in[k]$. We have 
$$\dg_{\cH(i)}(\alpha)=\sum\nolimits_{j\in [m]} a_{ij} h_j=nr_i,$$
and
$$\omega_\al(\cH(i))=\sum\nolimits_{j\in [m]} a_{ij}.$$ 
Therefore, if \eqref{longfatnec3} holds, 
$$\dg_{\cH(i)}(\al)-\omega_\al(\cH(i))=nr_i-\sum\nolimits_{j\in [m]} a_{ij}\geq n-1.$$
By Theorem \ref{mainthconfdet}, there exists  an $(\alpha,n)$-detachment $\cF$ of $\cH$ such that the following conditions hold.
\begin{enumerate} [label=(\alph*)]
\item For each $i\in [n]$ and each $j\in [k]$,
\begin{align*}
\dg_{\cF(j)} (\alpha_i)&\approx \frac{\dg_{\cH(j)}(\alpha)}{n}=\frac{nr_j}{n}=r_j;
\end{align*}
\item  For each $U\subseteq \{\alpha_1,\dots,\alpha_n\}$ with $|U|=h_j, j\in [m]$, 
\begin{align*}
\mult_\cF(U)&\approx  \frac{\mult_\cH(\alpha^{|U|})}{\binom{n}{|U|}}=\frac{\lambda_j \binom{n}{h_j}}{\binom{n}{h_j}}=\lambda_j;
\end{align*}
\item  For each $i\in [k]$, if \eqref{longfatnec3} holds, then $c(\cF(i))=1$. 
\end{enumerate}
By (b), $\cF\cong \Lambda K_n^H$, by (a) $\cF(i)$ is an $r_i$-factor, and by (c), for each $i$, if \eqref{longfatnec3} holds, then $\cF(i)$  is connected. 
\end{proof}
If we forbid edges of size 1, then condition \eqref{longfatnec3} will become simpler. 
\begin{corollary} \label{r1rkfacknhlamconcor}
$\Lambda K_n^H$ is $R$-factorable if and only if \eqref{longfatnec1} and \eqref{longfatnec2} are satisfied. Moreover, if $\min H\geq 2$, then for each $i\in [k]$, an $r_i$-factor is connected if and only if one of the following hold.
\begin{enumerate} [label=\textup{(\roman*)}]
\item $r_i=1, h_m=n,a_{im}=1, a_{ij}=0$ for $j\in [m-1]$,
\item $r_i\geq 2$.
\end{enumerate}
\end{corollary} 
\begin{proof}
Since $\min H\geq 2$, we have $2\leq h_1<h_2<\dots <h_m\leq n$. Fix $i\in [k]$.

Let us assume that \eqref{longfatnec3} holds. If $r_i\geq 2$, there is nothing to prove. If $r_i=1$, then by \eqref{longfatnec3}, $\sum\nolimits_{j\in[m]}  a_{ij}\leq 1$. Thus, there exists  an $\ell\in [m]$ such that $a_{i\ell}=1, a_{ij}=0$ for $j\neq \ell$.   Moreover, by \eqref{longfatnec2}, $h_\ell=\sum\nolimits_{j\in[m]} a_{ij} h_j=n$, and so $\ell=m$. 

Now, we show that if (i) or (ii) holds, then \eqref{longfatnec3} is satisfied. If (i) holds, then clearly   \eqref{longfatnec3} holds (It is also easy to see that if $h_m=n,a_{im}=1$, then the corresponding $r_i$-factor is connected). Let us assume that $r_i\geq 2$. By \eqref{longfatnec2},
$$
\frac{r_i n}{h_1}=\frac{1}{h_1} \sum\nolimits_{j\in[m]}  a_{ij} h_j=a_{i1}+\sum\nolimits_{j\in[m]} \frac{h_j}{h_1} a_{ij}>\sum\nolimits_{j\in[m]}  a_{ij}.
$$
Therefore, it suffices to show that $r_i n/h_1\leq n(r_i-1)+1$, or equivalently, $r_i\geq \frac{h_1(n-1)}{n(h_1-1)}$. Since $r_i\geq 2$, it is enough to show that $2n(h_1-1)\geq h_1(n-1)$, which is true because $h_1\geq 2$. 
\end{proof}

For two column vectors $Q=[q_1\dots q_k]^T$, $R=[r_1\dots r_k]^T$, we say that  $Q\leq R$ if $q_i\leq r_i$ for $i\in [k]$.  
A $(q,r)$-factor in a hypergraph $\cG$ is a spanning sub-hypergraph in which 
$$
q\leq \dg(v) \leq r \mbox { for each } v\in V(\cG).
$$
A \textit{$(Q,R)$-factorization} of $\cG$ is a partition $\{F_1,\ldots, F_k\}$ of $E(\cG)$ in which $F_i$ is a $(q_i,r_i)$-factor for $i\in  [k]$. 

The proof of the following result is very similar to that of Theorem \ref{r1rkfacknhlamcon}, and we shall skip it here.
\begin{theorem} \label{qrfacknhlamconqr}
$\Lambda K_n^H$ is $(Q,R)$-factorable if and only if
\begin{equation} \label{longfatnec1qr}
\sum\nolimits_{i\in [k]}q_i\leq \sum\nolimits_{i\in[m]}  \lambda_i\binom{n-1}{h_i-1}\leq \sum\nolimits_{i\in [k]}r_i,
\end{equation}
and there exists a non-negative integer matrix $A=[a_{ij}]_{k \times m}$ such that 
\begin{align} \label{longfatnec2qr}
    &&
   nQ\leq  AH\leq nR,
    &&
    \sum\nolimits_{i\in[k]}a_{ij}=\lambda_j\binom{n}{h_j} \mbox{ for } j\in [m].
    &&
  \end{align}
Moreover, for each $i\in [k]$, a $(q_i,r_i)$-factor that contains exactly $a_{ij}$ edges of size $h_j$ for each $j\in [s]$, is connected if and only if
\begin{equation} \label{longfatnec3qr}
\sum\nolimits_{j\in[m]}  a_{ij}(h_j-1)\geq n-1.
\end{equation}
\end{theorem}
An \textit{almost $k$-factor} of a hypergraph is $(k-1,k)$-factor. Theorem \ref{qrfacknhlamconqr} can be used to construct connected almost factorizations of $\Lambda K_n^H$. 

\section*{Acknowledgement}
We are deeply grateful to both referees for their  careful reading of the manuscript and for their numerous insightful and constructive suggestions, which have greatly improved the paper. 
 The   author's research is  partially supported by  NSA Grant H98230-16-1-0304.

\bibliographystyle{plain}  
   \bibliography{connfairbib}
\end{document}